\documentclass{amsart}
\usepackage{graphicx}

\newcommand{\tens}{\otimes}
\newcommand{\join}{\vee}
\newcommand{\meet}{\wedge}
\newcommand{\derc}[1]{C^{(#1)}}
\newcommand{\derm}[1]{M^{(#1)}(C)}
\newcommand{\cut}[1]{}
\newcommand{\ext}[1]{$#1$-extensor}
\newcommand{\exts}[1]{$#1$-extensors}
\newcommand{\drop}[1]{(------ {\bf (text omitted) ------}}
\newcommand{\leave}[1]{}
\newcommand{\gcr}{Gian Carlo Rota}
\newcommand{\gc}{Gian Carlo}
\newcommand{\cop}[3]{$\partial_{#1,#2}\;#3$}
\newcommand{\copm}[3]{\partial_{#1,#2}\;#3}
\newcommand{\quo}[1]{{\it``#1''}}
\newcommand{\ra}[1]{\rho(#1)}
\newcommand{\pr}{\circ}
\newcommand{\dott}[1]{\dot{{\phantom |}#1{\phantom |}}}
\newcommand{\sgn}{\sigma}
\newcommand{\col}[1]{#1} 
\newcommand{\colAlt}[2]{#1} 
\begin{document}
\title[An Algebra of Pieces of Space]{An Algebra of Pieces of Space ---\\ 
Hermann Grassmann to Gian Carlo Rota }
\author{Henry Crapo}
\email{crapo@ehess.fr}
\address{Centre de Recherche ``Les Moutons matheux'', 
34520 La Vacquerie, France.}
\begin{abstract}
We sketch the outlines of \gcr's interaction with the
ideas that Hermann Grassmann developed in
his {\it Ausdehnungslehre}\cite{grass44,grass62}
 of 1844 and 1862, as adapted 
and explained by Giuseppe Peano in 1888. This leads us past what \gc\
variously called {\it Grassmann-Cayley algebra} and {\it Peano spaces}
to the {\it Whitney algebra} of a matroid, and finally to a resolution
of the question ``What, really, was Grassmann's {\it regressive product}?''.
This final question is the subject of ongoing joint work with
Andrea Brini, Francesco Regonati, and William Schmitt. 
\end{abstract}
\maketitle

\leave{Alternate title: From Grassmann to Rota and beyond, a voyage in time and spaces.}
\section{Almost ten years later}
We are gathered today in order to renew and deepen our recollection
of the ways in which our paths intersected that of \gcr.
We do this in poignant sadness, but with a bitter-sweet touch: we are 
pleased to have this opportunity to meet and to discuss his life and work, since
we know how \gc\ transformed us through his friendship and his love of
mathematics.

We will deal only with the most elementary of geometric questions; how to 
represent pieces of space of various dimensions, in their relation to one another.
It's a simple story, but one that extends over a period of some 160 years.
We'll start and finish with Hermann Grassmann's project, but the trail will lead
us by Giuseppe Peano, Hassler Whitney, to \gcr\ and his colleagues.

Before I start, let me pause for a moment to recall a late afternoon at the
Accademia Nazionale dei Lincei, in 1973, on the eve of another talk 
I was petrified to give, when
\gc\ decided to teach me how to {\it talk}, so I wouldn't make a fool
of myself the following day. The procedure was for me to start my talk, 
with an audience of one, and
he would interrupt whenever there was a problem. We were in that otherwise empty
conference hall for over two hours, and I never got past my first paragraph.
It was terrifying, but it at least got me through the first battle with my fears
and apprehensions, disguised as they usually are by timidity, self-effacement,
and other forms of apologetic behavior.

\section{Synthetic Projective Geometry}

Grassmann's plan was to develop a purely formal algebra to model 
natural (synthetic) operations on geometric objects: 
{\it flat}, or {\it linear} pieces of space of all possible dimensions.
His approach was to be {\it synthetic}, so that the symbols in his
algebra would denote geometric objects themselves, not just numbers
(typically, coordinates) that could be derived from those objects
by measurement. His was not to be an algebra of numerical quantities,
but an algebra of pieces of space.

In the {\it analytic} approach, so typical in the
teaching of Euclidean geometry, we are encouraged to assign 
``unknown'' variables to the coordinates of variable points,
to express our hypotheses as equations in those coordinates,
and to derive equations that will express our desired
conclusions.

The main advantage of a synthetic approach is that the logic of
geometric thought and the logic of algebraic manipulations may
conceivably remain parallel, and may continue to cast light upon 
one another. Grassmann expressed this clearly
in his introduction to the {\it Ausdehnungslehre}\cite{grass44,grass44T}:

Grassmann (1844): \quo{Each step from one formula to another appears at once
as just the symbolic expression of a parallel act of abstract reasoning. 
The methods formerly used require the introduction of arbitrary
coordinates that have nothing to do with the problem and completely
obscure the basic idea, leaving the calculation as a mechanical 
generation of formulas, inaccessible and thus deadening to the intellect. 
Here, however, where the idea is no longer strangely obscured but
radiates through the formulas in complete clarity, the intellect
grasps the progressive development of the idea with each formal
development of the mathematics.}

In our contemporary setting, a synthetic approach to geometry
yields additional benefits. At the completion of a synthetic calculation,
there is no need to {\it climb back up} from scalars (real numbers, 
necessarily subject to round-off errors, often rendered totally
useless by division by zero) or from drawings, fraught with
their own approximations of incidence, to statements of geometric
incidence. In the synthetic approach, one even receives precise warnings
as to particular positions of degeneracy. The synthetic approach
is thus tailor-made for machine computation.

\gc\ was a stalwart proponent of the synthetic approach to geometry
during the decade of the 1960's, when he studied the combinatorics
of ordered sets and lattices, and in particular, matroid theory.
But this attitude did not withstand his encounter with invariant theory, beginning
with his lectures on the invariant theory of the symmetric group at the A.M.S.
summer school at Bowdoin College in 1971.

As \gc\ later said, with his admirable fluency of expression in
his adopted tongue,

\quo{Synthetic projective geometry in the plane held great sway
between 1850 and 1940. It is an instance of a theory whose beauty was 
largely in the eyes of its beholders. Numerous expositions were written 
of this theory by English and Italian mathematicians (the definitive one being
the one given by the American mathematicians Veblen and Young). These
expositions vied with one another in elegance of presentation and in
cleverness of proof; the subject became required by universities
in several countries. In retrospect, one wonders what all the fuss was about.}

\quo{Nowadays, synthetic geometry is largely cultivated by historians, and
the average mathematician ignores the main results of this once flourishing
branch of mathematics. The claim that has been raised by defenders of synthetic
geometry, that synthetic proofs are more beautiful than analytic proofs,
is demonstrably false. Even in the nineteenth century, invariant-theoretic
techniques were available that could have provided elegant, coordinate-free 
analytic proofs of geometric facts without resorting to the gymnastics
of synthetic reasoning and without having to stoop to using coordinates.}

Once one adopts an invariant-theoretic approach, much attention must be paid
to the reduction of expressions to standard form, where one can recognize
whether a polynomial in determinants is equal to zero. The process is called
{\it straightening}, and it was mainly to the algorithmic process of straightening
in a succession of algebraic contexts that \gc\ devoted his creative 
talents during three decades. We filled pages with calculations such as
the following, for the bracket algebra of six points in a projective plane:
$$
\begin{matrix}
b&c&d\\
a&e&f
\end{matrix}
\;\;-\;\;
\begin{matrix}
a&c&d\\
b&e&f
\end{matrix}
\;+\;
\begin{matrix}
a&b&d\\
c&e&f
\end{matrix}
\;-\;
\begin{matrix}
a&b&c\\
d&e&f
\end{matrix}
\;\;=\;\;0
$$
the {\it straightening} of the two-rowed tableau on the left. 
This expression we would write
in {\it dotted} form
$$
\begin{matrix}
\dott b&\dott c&\dott d\\
\dott a&e&f
\end{matrix}
$$
that would be {\it expanded} to the above expression by summing, with
alternating sign, over all permutations of the dotted letters. The basic principle is that {\it dotting a dependent set of points yields zero.}

To make a long and fascinating story short,  \gc\ finally settled upon 
a most satisfactory compromise,
a formal super-algebraic calculus, developed with Joel Stein, Andrea Brini,
Marilena Barnabei \cite{bbr,bb,bht,brt}
 and a generation of graduate students, 
that managed to stay reasonably close
to its synthetic geometric roots. His brilliant students Rosa Huang and
Wendy Chan \cite{chan3-6,crs, bht, hrs}
carried the torch of synthetic reasoning across this new territory.
They rendered feasible a unification of super-algebraic and synthetic geometry, 
but, as we soon realize, the process is far from complete.

First, we should take a closer look at Grassmann's program for synthetic geometry.

\section{Hermann Grassmann's algebra}
Grassmann emphasizes that he is building an {\it abstract} theory that can then
be {\it applied to real physical space}. He starts not with geometric axioms,
but simply with the notion of a skew-symmetric product of letters, which
is assumed to be distributive over addition and modifiable by scalar
multiplication. Thus if $a$ and $b$ are letters, and $A$ is any product of 
letters, then $Aab + Aba = 0$. This is the skew-symmetry. It follows
that any product with a repeated letter is equal to zero.

He also develops a notion of {\it dependency}, saying that a letter $e$
is {\it dependent} upon a set $\{a,b,c,d\}$ 
if and only if
there are scalars $\alpha,\beta,\gamma,\delta$ such that
$$\alpha\, a + \beta\, b + \gamma\, c+ \delta\, d = e.$$
Grassmann realizes that such an expression is possible if and only if
the point $e$ lies in the {\it common system}, or projective subspace
spanned by the points $a,b,c,d$. He
uses an axiom of distributivity of product over sums to prove that the product
of letters forming a dependent set is equal to zero.
With $a,b,c,d,e$ as above:
$$abcde = 
 \alpha\, abcda + \beta\, abcdb + \gamma\, abcdc + \delta\, abcdd = 0$$
the terms on the right being zero as products because each has a 
repeated letter.

As far as I can see, 
he establishes no formal axiomatization of the relation of 
linear dependence, and in particular, no statement of the exchange property.
For that, we must wait until 1935, and the matroid theory of Whitney,
MacLane and Birkhoff.

The application to geometry is proposed via {\it an interpretation}
of this abstract algebra. The individual letters may be understood as
{\it points},

{\it The center of gravity of several points can be interpreted 
as their sum, the displacement between two points as their product, 
the surface area lying between three points as their product, and 
the region (a pyramid) between four points as} {their} {\it product.}

\cut {\it For a long time it had been evident to me that geometry can in 
no way be regarded as a branch of mathematics like arithmetic or 
combination theory; instead, geometry relates to something already 
given in nature, namely
space. I had also realized that there must be a branch of mathematics that
yields in a purely abstract way laws similar to those that in geometry seem 
bound to space. The possibility of constructing such a purely abstract 
branch of mathematics was provided by the new analysis; indeed this analysis, 
developed without the assumption of any principles established outside
its own domain, and proceeding entirely by abstraction, 
was itself this science.}

Grassmann is delighted to find that, in contrast to earlier formalizations of 
geometry, there need be no a priori maximum to the rank of the 
overall space being described:

{\it The origins of this science are as immediate as those of arithmetic;
and in regard to content, the limitation to three dimensions is absent.
Only thus do the laws come to light in their full clarity, and their essential
interrelationships are revealed.}

\begin{figure}[h] 
     \centering
     \includegraphics{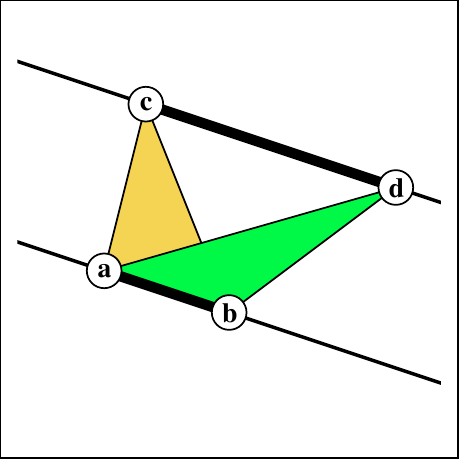}
     \caption{$abc=abd$ {\it iff} line $ab$ parallel to line $cd$.}
     \label{Fi:area}
     \end{figure}

Giuseppe Peano, in rewriting Grassmann's theory, chose to assume a few basic 
principles of {\it comparison of signed areas and volumes}, 
and to base all proofs
on those principles. For instance, given a line segment $ab$ 
between distinct points $a,b$ in the plane,
and points $c,d$ off that line, then the signed areas $abc$ and $abd$ will
be equal (equal in magnitude, and equally oriented CW or CCW) if and only
if $c$ and $d$ lie on a line {\it parallel} to $ab$. See Figure~\ref{Fi:area}.

\begin{figure}[h] 
     \centering
     \includegraphics{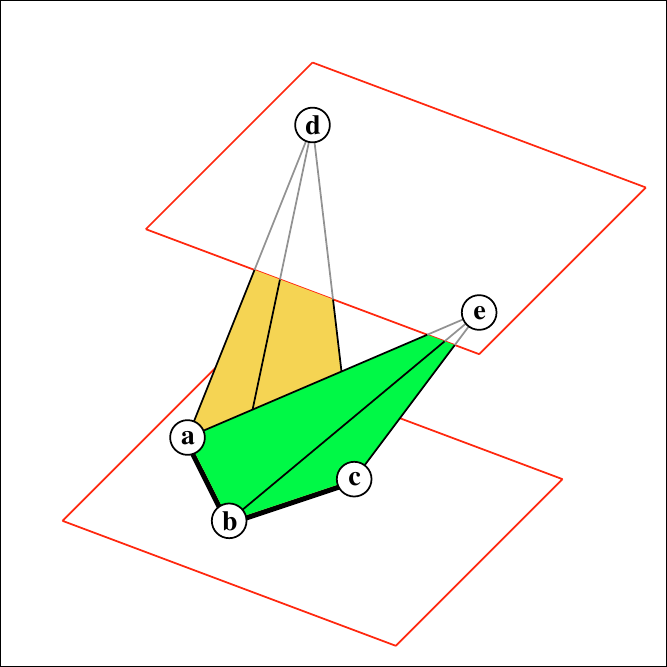}
     \caption{$abcd = abce$ {\it iff} line $de$ parallel to plane $abc$.}
     \label{Fi:vol}
     \end{figure}

The corresponding statement
for three-dimensional space is that signed volumes $abcd$ and $abce$,
for non-collinear points $a,b,c$, are
equal (equal volume, and with the same chirality, or handedness) 
 if and only if $d$ and $e$ 
lie on a line {\it parallel} to the plane $abc$. See Figure~\ref{Fi:vol}. 
Since Peano restricts his attention
to three-dimensional space, this principle is his main axiom, with a notion
of parallelism taken to be understood. This means that even the simplest
geometric properties are ultimately rephrased as equations among measured
volumes. For instance, Peano wishes to show that three points $a,b,c$ are
{\it collinear} if and only if the linear combination  
$bc - ac + ab$ of products of points is equal to zero.
He shows that for every pair $p,q$ of points, if the points $a,b,c$
are in that order on a line, the tetrahedron $acpq$ is the disjoint
union of the tetrahedra $abpq$ and $bcpq$, so their volumes add:
$$
abpq + bcpq = acpq.
$$
A further argument about symmetry shows that $ab + bc = ac$ holds independent
of the order of $a,b,c$ on the line. The statement $abpq + bcpq = acpq$,
quantified over all choices of
$p$ and $q$, is Peano's {\it definition} of the equality $ab + bc = ac$.
So his proof is complete.
Perhaps we can agree that this is putting the cart before the horse.
(Grassmann took expressions of the form $ab + bc = ac$, for three
collinear points $a,b,c$, to be axiomatic in his algebra.)

I mention this strange feature of Peano's version because it became
something of an {\it id\'ee fixe} for \gc's work on Peano spaces.
It gave rise to the technique of {\it filling brackets} in order to verify
equations involving flats of low rank, a technique that
unnecessarily relies on information concerning the rank of the overall space.

\section{Extensors and Vectors}
We are all familiar with the formation of linear combinations of 
points in space, having studied linear algebra and having learned
to replace points by their {\it position vectors}, which can then
be multiplied by scalars and added.  The origin serves as reference
point, and that is all we need.

Addition of points is not well-defined in real projective geometry, because
although points in a space of rank $n$ (projective dimension $n-1$)
may be represented as $n$-tuples of real numbers, the $n$-tuples are
only determined up to an overall non-zero scalar multiple, and addition
of these vectors will not produce a well-defined result. The usual
approach is to consider {\it weighted points}, consisting of a
position, given by {\it standard homogeneous coordinates}, of the form
$(a_1,\dots,a_{n-1},1)$, and a {\it weight} $\mu$, to form a 
{\it point $(\mu a_1,\dots,\mu a_{n-1},\mu)$ of weight $\mu$}.
This is what worked for M\"obius in his {\it barycentric calculus}.
And it is the crucial step used by Peano to clarify the presentation
of Grassmann's algebra.

\begin{figure}[h] 
     \centering
     \includegraphics{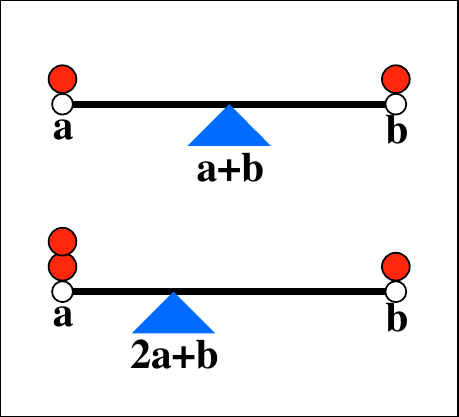}
     \caption{Addition of weighted points.}
     \label{Fi:weights}
     \end{figure}

This amounts to fixing a choice of {\it hyperplane at infinity} with
equation $x_n = 0$. The {\it finite points} are represented as 
above, with weight 1. A linear combination $\lambda a + \mu c$, for
scalars $\lambda$ and $\mu$ positive, 
is a point $b$ situated between $a$ and $c$, such that the ratio
of the distance  $a\to b$ to the distance $b \to c$ is in the
inverse proportion $\mu/\lambda$, and the resulting weighted
point has weight equal to 
$\lambda+\mu$, as illustrated in Figure~\ref{Fi:weights}. 

In particular, for points $a$ and $b$ of weight 1,
$a+b$ is located the midpoint of the
interval $ab$, and has weight $2$, while $2a+b$ is located twice
as far from $b$ as from $a$, and has weight $3$.

Both Grassmann and Peano are careful to distinguish between products $ab$ of
points, which have come to be called {\it \exts{2}}, and differences 
$b-a$ of points, which Grassmann calls {\it displacements} and Peano calls
{\it $1$-vectors}. The distinction between such types of objects is easily
explained in terms of modern notation, in homogeneous coordinates. 

\begin{figure}[h] 
     \centering
     \includegraphics{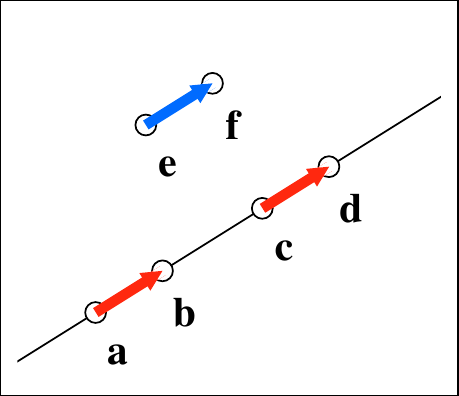}
     \caption{$ab=cd$, while $b-a=d-c=f-e$ is a point at infinity.}
     \label{Fi:twoExt}
     \end{figure}

In Figure~\ref{Fi:twoExt} we show two equal \exts{2}, $ab=cd$, in red, and their difference
vector, $v=b-a=d-c=f-e$ in blue, which represents a projective point, not a line, namely, the
point at infinity on the line $ab$, with weight equal to the length from $a$ to $b$ and
sign indicating the orientation form $a$ to $b$. Check that $ab=av=bv$, so multiplication
of a point $a$ on the right by a vector $v$ creates a line segment of length and direction $v$
starting at $a$.

To avoid
cumbersome notation, we will henceforth follow Peano's example, and 
give all examples with reference
to a real projective $3$-space,  rank 4.

The homogeneous coordinates of any product of $k\le 4$ weighted points (called a {\it \ext{k}}), are 
the $k\times k$ minors of the matrix whose rows are the coordinate vectors of the
$k$ points in question. So for any four weighted points $a,b,c,d$ in a space of
rank $4$,

$$
a = \begin{pmatrix}
\bf1&\bf2&\bf3&\bf4\\
a_1&a_2&a_3&a_4
\end{pmatrix}
$$
$$
ab = \begin{pmatrix}
\bf12&\bf13&\bf23&\bf14&\bf24&\bf34\\
\begin{vmatrix}
a_1&a_2\\
b_1&b_2
\end{vmatrix}&
\begin{vmatrix}
a_1&a_3\\
b_1&b_3
\end{vmatrix}&
\begin{vmatrix}
a_2&a_2\\
b_2&b_3
\end{vmatrix}&
\begin{vmatrix}
a_1&a_4\\
b_1&b_4
\end{vmatrix}&
\begin{vmatrix}
a_2&a_4\\
b_2&b_4
\end{vmatrix}&
\begin{vmatrix}
a_3&a_4\\
b_3&b_4
\end{vmatrix}
\end{pmatrix}
$$

$$
abc = 
\begin{pmatrix}
\bf123&\bf124&\bf134&\bf234\\
\begin{vmatrix}
a_1&a_2&a_3\\
b_1&b_2&b_3\\
c_1&c_2&c_3
\end{vmatrix}&
\begin{vmatrix}
a_1&a_2&a_4\\
b_1&b_2&b_4\\
c_1&c_2&c_4
\end{vmatrix}&
\begin{vmatrix}
a_1&a_3&a_4\\
b_1&b_3&b_4\\
c_1&c_3&c_4
\end{vmatrix}&
\begin{vmatrix}
a_2&a_3&a_4\\
b_2&b_3&b_4\\
c_2&c_3&c_4
\end{vmatrix}
\end{pmatrix}
$$
$$
abcd =
\begin{pmatrix}
\bf1234\\
\begin{vmatrix}
a_1&a_2&a_3&a_4\\
b_1&b_2&b_3&b_4\\
c_1&c_2&c_3&c_4\\
d_1&d_2&d_3&d_4
\end{vmatrix}
\end{pmatrix}
$$

If $a,b,c,d$ are points of weight 1, the \exts{2}\ $ab$ and $cd$ are equal if and only
if the line segments from $a$ to $b$ and from $c$ to $d$ are {\it collinear}, of equal length,
and similarly oriented. More generally, Grassmann showed that two \exts{k}\ $a\dots b$
and $c\dots d$ differ only by a non-zero scalar multiple, $a\dots b = \sigma\; a\dots d$, 
if and only if the sets of points obtainable as linear combinations of $a,\dots, b$
and those from $c,\dots,d$ form what we would these days call {\it the same projective 
subspace}. Such a subspace, considered as a set of projective points, we call a 
{\it projective flat}.

\begin{figure}[h] 
     \centering
     \includegraphics{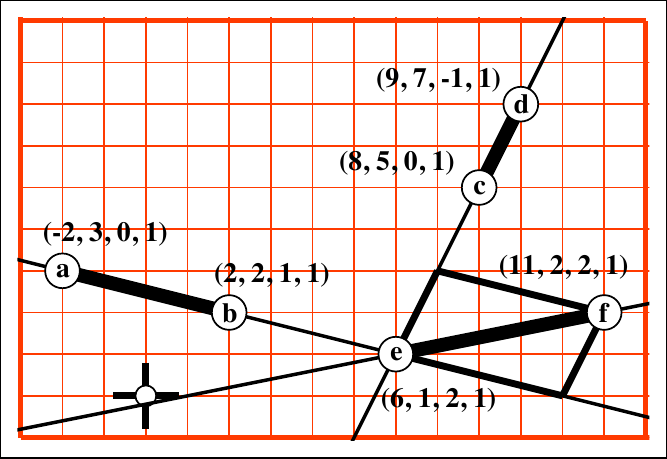}
     \caption{$ab+cd=ef$ for these six points of weight $1$.}
     \label{Fi:add}
     \end{figure}

Coplanar \exts{2}\ add the way coplanar forces do in physical systems. 
Say you are forming the sum $ab+cd$ as in Figure~\ref{Fi:add}. You slide 
the line segments representing the forces $ab$ and $cd$
along their lines of action until  the ends $a$ and $c$ coincide at the point $e$ of incidence
of those two lines. The sum is then represented as the diagonal line segment of the parallelogram
they generate.

Let's carry out the explicit extensor calculations, 
so you know I'm not bluffing: let $a,b,c,d$ be four points coplanar in $3$-space
(rank $4$),
$$
\begin{matrix}
&\bf1&\bf2&\bf3&\bf4\\
\bf a&-2&3&0&1\\
\bf b&2&2&1&1\\
\bf c&8&5&0&1\\
\bf d&9&7&-1&1
\end{matrix}
$$
That the four points are coplanar is clear from the fact that
the four vectors are dependent: $-a+2b-3c+2d=0$.
The \exts{2}\ $ab$ and $cd$ are, with their sum:
$$
\begin{matrix}
&\bf12&\bf13&\bf23&\bf14&\bf24&\bf34\\
\bf ab&-10&-2&3&-4&1&-1\\
\bf cd&11&-8&-5&-1&-2&1\\
\bf ab+cd&1&-10&-2&-5&-1&0
\end{matrix}
$$
The point $e$ of intersection of lines $ab$ and $cd$, together with the 
point $f$ situated at the end of the diagonal of the parallelogram formed
by the translates of the two line segments to $e$, have homogeneous 
coordinates
$$
\begin{matrix}
&\bf1&\bf2&\bf3&\bf4\\
\bf e&6&1&2&1\\
\bf f&11&2&2&1
\end{matrix}
$$
and exterior product
$$
\begin{matrix}
&\bf12&\bf13&\bf23&\bf14&\bf24&\bf34\\
\bf ef&1&-10&-2&-5&-1&0
\end{matrix}
$$
equal to the sum $ab+cd$.

\begin{figure}[h] 
     \centering
     \includegraphics{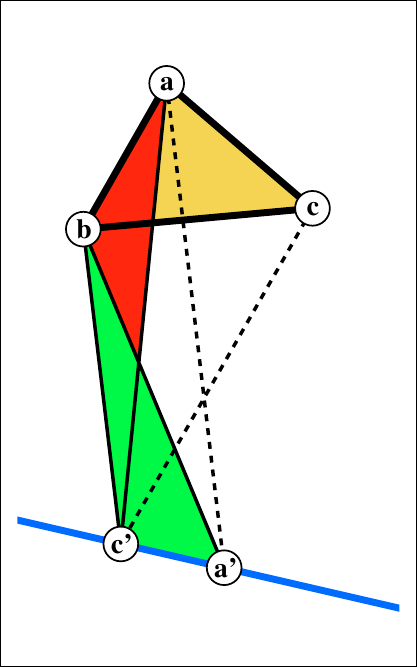}
     \caption{The \exts{3} $abc$ and $a'bc'$ are equal.}
     \label{Fi:pla}
     \end{figure}

Cospatial planes add in a similar fashion. For any \ext{3}\ $abc$ spanning a
projective plane $Q$ and for any 
\ext{2}\ $pq$ in the plane $Q$, there are points $r$ in $Q$ such that $abc=pqr$.
The required procedure is illustrated in Figure~\ref{Fi:pla}. We slide the point
$c$ parallel to the line $ab$ until it reaches the line $pq$, shown in blue, 
at $c'$. Then slide
the point $a$ parallel to the line $bc'$ until it reaches the line $pq$ at $a'$.
The oriented plane areas $abc$, $abc'$, and $a'bc'$ are all equal, and the final
triangle has an edge on the line $pq$.

\begin{figure}[h] 
     \centering
     \includegraphics{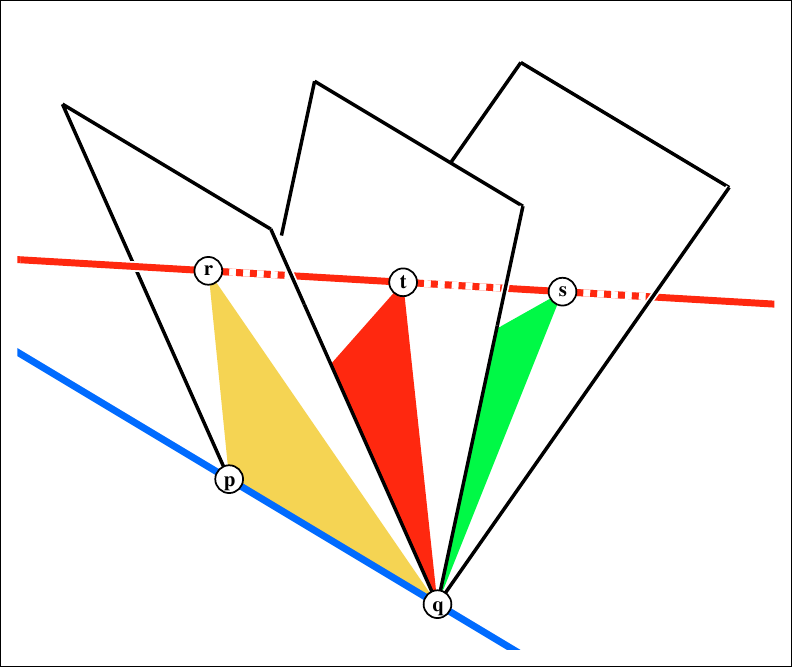}
     \caption{A sum of cospatial \exts{3}, $pqr+pqs=pqt$.}
     \label{Fi:addPla}
     \end{figure}

So, given any \exts{3}\ $abc,def$, and for any pair $p,q$ of points
on the line of intersection of the planes $abc$ and $def$, there exist points 
$r,s$ such that $abc=pqr, def=pqs$, so $abc+def$ can be expressed in the form
$pqr+pqs=pq(r+s)$, and the problem of adding planes in $3$-space is reduced
to the problem of adding points on a line. The result is shown in 
Figure~\ref{Fi:addPla}, where $t=r+s$.

In dimensions $\ge 3$, it is necessary to {\it label} the individual coordinates
with the set of columns used to calculate that minor. This practice 
becomes even more systematic in the subsequent
super-algebraic approach of Rota, Stein, Brini and colleagues,
where a coordinate $abc_{ijk}$ will be 
denoted $(abc|ijk)$, with {\it negative letters} $a,b,c$ to denote vectors 
and {\it negative
places} $i,j,k$ to denote coordinate positions, and the minor is 
calculated by the {\it Laplace expansion}:
$$\begin{matrix}
(abc|ijk)&=& (a|i)(b|j)(c|k) - (a|i)(b|k)(c|j) + \dots - (a|k)(b|j)(c|i),\\
&=&(a|i)(b|j)(c|k) - (a|i)(c|j)(b|k) +\dots -(c|i)(b|j)(a|k)
\end{matrix}$$
as the sum of products of individual {\it letter-place} elements.
We have chosen to list the coordinates of a line segment $ab$ in $3$-space in the order
$12, 13, 23, 14, 24, 34$ because that makes the negative of the difference vector $b-a$
visible in the last three coordinate places $14,24,34$, and the {\it moments} about the 
$3^{rd}$, $2^{nd}$, and $1^{st}$ coordinate axes,
respectively, visible in the first three coordinate places. 

A set of $k$ \exts{1}\ $a,b,\dots d$ has a non-zero product if and only if the $k$
points are independent, and thus span a projective subspace of rank $k$ (dimension $k-1$).
This integer $k$ is called the {\it step} of the extensor.

We have seen that differences of points $a,b$ of weight $1$ are vectors, which means
simply that they are \exts{1}\ $e$ {\it at infinity}, with coordinate $e_4 = 0$.

In Grassmann's theory, there exist $k$-vectors of all steps $k$ for which
$k$ is less than the rank of the entire space. 
In terms of standard homogeneous coordinates in a space of rank $4$, they are
those extensors for which all coordinates are zero whose labels involve the place $4$.
$k$-vectors are also definable as ``boundaries'' of \exts{(k+1)}, or as products
of $1$-vectors, as we shall see.

Each \ext{k}\ has an associated $(k-1)$-vector, which I shall refer to as its {\it boundary},
$$\begin{matrix}
\partial\, ab = b-a\\
\partial\, abc = bc - ac + ab\\
\partial\, abcd = bcd - acd + abd - abc
\end{matrix}
$$
If the points $a,b,c,d$ are of weight 1, then the $4^{th}$ coordinate of $\partial\, ab$,
the $14, 24, 34$ coordinates of $\partial\, abc$, and the $124, 134, 234$ coordinates of $\partial\, abcd$
are all equal to zero. Such extensors are, as elements of our affine version of projective space, 
pieces of space in the hyperplane at infinity. 
As algebraic objects they are vectors, so $\partial\, ab = \partial\, cd$,
or $b-a = d-c$, for points $a,b,c,d$ of weight 1, if and only if the line segments from $a$ to $b$
and from $c$ to $d$ are parallel, of equal length, and similarly oriented. That is, they are equal as
difference vectors of position. 

The subspace of $k$-vectors consists exactly of those $k$-tensors obtained by taking 
boundaries of \exts{k+1}. The $k$-vectors are also
expressible as  exterior products
of $k$ $1$-vectors, since 
$$\begin{matrix}
\partial\, ab = (b-a)\\
\partial\, abc = (b-a)(c-a),\\
 \partial\, abcd = (b-a)(c-a)(d-a),\\
 \dots.
\end{matrix}
$$
Remarkably, $p\,\partial\, abc = abc$ for any point $p$ in the plane spanned by $a,b,c$.
In particular, $a\,\partial\, abc = b\,\partial\, abc = c\,\partial\, abc$.

So $\partial\, abc = bc - ac + ab$ is a $2$-vector, and it can
only have non-zero coordinates with labels $12, 13, 23$. Geometrically, it can
be considered to be a directed line segment in the line at infinity in the plane $abc$.
It is also equal to a {\it couple} in the plane of the points $a,b,c$, that is, the 
sum of two \exts2\ that are parallel to one another, of equal length and opposite
orientation.
Couples of forces occur in statics,
and cause rotation when applied to a rigid body.

\begin{figure}[h] 
     \centering
     \includegraphics{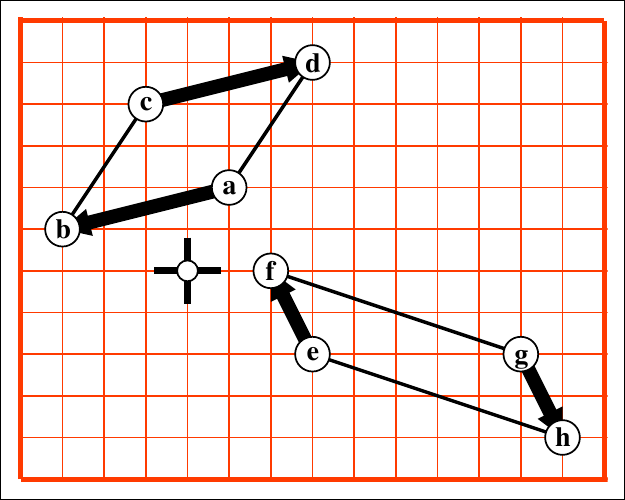}
     \caption{The couples $ab+cd$ and $ef+gh$ are equal.}
     \label{Fi:cou}
     \end{figure}

In Figure~\ref{Fi:cou}
 we show two couples that are equal, though expressed
as sums of quite different \exts{2}. The coordinate expressions are as follows.
$$
\begin{matrix}
&\bf1&\bf2&\bf3\\
\bf a&1&2&1\\
\bf b&-3&1&1\\
\bf c&-1&4&1\\
\bf d&3&5&1
\end{matrix}\quad\quad\quad\quad
\begin{matrix}
&\bf1&\bf2&\bf3\\
\bf e&3&-2&1\\
\bf f&2&0&1\\
\bf g&8&-2&1\\
\bf h&9&-4&1
\end{matrix}
$$
$$
\begin{matrix}
ab+cd=
\begin{pmatrix}
\bf12&\bf13&\bf23\\
7&4&1
\end{pmatrix}
+
\begin{pmatrix}
\bf12&\bf13&\bf23\\
-17&-4&-1
\end{pmatrix}
=
\begin{pmatrix}
\bf12&\bf13&\bf23\\
-10&0&0
\end{pmatrix}\\
ef+gh=
\begin{pmatrix}
\bf12&\bf13&\bf23\\
4&1&-2
\end{pmatrix}
+
\begin{pmatrix}
\bf12&\bf13&\bf23\\
-14&-1&2
\end{pmatrix}
=
\begin{pmatrix}
\bf12&\bf13&\bf23\\
-10&0&0
\end{pmatrix}
\end{matrix}
$$

$$
\begin{matrix}
ab+cd = (7,4,1) + (-17,-4,-1) = (-10,0,0)\\
ef+gh = (4,1,-2) + (-14,-1,2) = (-10,0,0)
\end{matrix}
$$

\section{Reduced Forms}

Grassmann did much more than create a new algebra of {\it pieces of space}. First he showed 
that every sum of \exts{1}\ (weighted points), with sum of weights not equal to zero,
is equal to a single weighted point (with weight equal to the sum of the weights
of the individual points). If the sum of the weights is zero, the resulting \ext{1}\ is
a $1$-vector (which may itself be zero). He then went on to show that every sum of coplanar 
\exts{2}\ is equal to a single \ext{2}, or to the sum of two \exts{2}\
on parallel lines, with equal length and opposite orientation, forming a couple,
or is simply zero.

The situation for linear combinations of \exts{2}\ in $3$-space is a bit more complicated.
The sum of \exts{2}\ that are not coplanar is not expressible as a product of points,
and so is not itself an extensor.
It is simply an antisymmetric tensor of step $2$.
Grassmann showed that such a linear combination of non-coplanar extensors
can be reduced to the sum of a \ext{2}\ and a $2$-vector, or couple. 

Not bad at all, for the mid 19$^{th}$ century.

At the beginning of
the next millenium, in 1900, Sir Robert Ball, in his {\it Theory of Screws} will use a bit
of Euclidean geometry to show that a sum of \exts{2}\ in $3$-space can be expressed
as the sum of a force along a line plus a moment in a plane perpendicular to that line.
He called these general antisymmetric tensors {\it screws}.
Such a combination of forces, also called a {\it wrench}, when applied to a rigid body produces a {\it screw motion}.

Much of the study of screws, with applications to statics, is to be found at the end of
Chapter 2 in Grassmann. He discusses coordinate notation, change of basis, and even
shows that an anti-symmetric tensor $S$ of step $2$  (a screw) is an extensor if and only if the 
exterior product $SS$ is equal to zero. This is the first and most basic invariant
property of anti-symmetric tensors.

Any linear combination of \exts{3}\ in rank $4$ ($3$-space) is equal to a single \ext{3}.
This is because we are getting close to the rank of the entire space. The simplicity of
calculations with linear combinations of points is carried over by duality to calculations
with linear combinations of copoints, here, with planes.

The extreme case of this duality becomes visible in that \exts{k}\ in $k$-space add and
multiply just like scalars. For this reason they are called {\it pseudo-scalars}.

\section{Grassmann-Cayley Algebra, Peano Spaces}
\gc\ chose to convert
pseudo-scalars to ordinary scalars by taking a determinant, or {\it bracket} of the product\cite{drs}.
These brackets provide the scalar coefficients of any invariant linear combination. Thus,
for any three points $a,b,c$ on a projective line, there is a linear relation, 
unique up to an overall scalar multiple
$$
[ab] c - [ac] b + [bc] a = 0
$$
and for any four points $a,b,c,d$ on a projective plane
$$
[abc] d - [abd] c + [acd] b - [bcd] a = 0
$$
So for any four coplanar points $a,b,c,d$, 
$$[abc] d - [abd] c = [bcd] a - [acd] b,$$
and we have two distinct expressions
for a projective point that is clearly on both of the lines $cd$ and $ab$.

\begin{figure}[h] 
     \centering
     \includegraphics{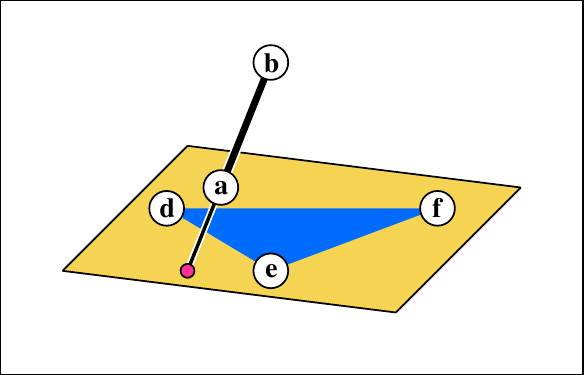}
     \caption{The meet of  $ab$ with $cde$.}
     \label{Fi:meet}
     \end{figure}

This is the key point in the development of the Grassmann-Cayley algebra
\cite{drs,grs, whiteCGA},
as introduced by \gc\ with his coauthors Peter Doubilet and Joel Stein.
If $A$ is an \ext{r}\ and $B$ is a \ext{s}\ in a space of overall rank $n$,
then the {\it meet} of $A$ and $B$ is defined by the two equivalent formulae
$$
A\meet B= \sum_{(A)_{r-k,k}}[A_{(1)}B]A_{(2)}= \sum_{(B)_{k,s-k}}[AB_{(2)}]B_{(1)}
$$
This is the {\it Sweedler notation} (see below, Section~\ref{Se:wa}) 
from Hopf algebra, where, for instance,
in the projective $3$-space, rank $4$, a line $ab$ and a plane $cde$
as in Figure~\ref{Fi:meet} will have meet
$$
ab\meet cde \;\;=\;\; [acde] b - [bcde] a \;\;=\;\; [deab] c - [ceab] d + [cdab] e.
$$
This meet is equal to zero unless $ab$ and $cde$ together span the whole space,
so the line meets the plane in a single point.
You can check that the second equality checks with the generic relation
$$
[abcd] e - [abce] d + [abde] c - [acde] b + [bcde] a \;\;=\;\; 0.
$$

The serious part of Grassmann-Cayley algebra starts when you try to use these
simple relations to detect properties of geometric configuration in {\it special
position}, or to classify possible interrelations between subspaces, that is,
when you start work in {\it invariant
theory}. 

\gc\ and his colleagues made considerable progress in expressing the better-known
invariants in terms of this new double algebra, and in reproving a certain number of
theorems of projective geometry. Matters got a bit complicated, however, when they got
to what they called the {\it alternating laws}, which ``alternated'' the operations 
of join (exterior product) and meet in a single equation. This required maintaining
a strict accounting of the ranks of joins and meets, in order to avoid unexpected
zeros along the way. They employed Peano's old technique of {\it filling brackets}
whenever the situation got delicate. It took some real gymnastics to construct
propositions of general validity when mixing joins and meets.

\section{Whitney algebra}\label{Se:wa}
The idea of the Whitney algebra of a matroid starts with the simple
observation that the relation
$$
a\tens bc - b\tens ac + c\tens ab = 0
$$
among tensor products of extensors holds if and only if the three
points $a,b,c$ are collinear, and this {\it in a projective space of arbitrary
rank.}

We need the notion of a {\it coproduct slice} of a word. For any non-negative
integers $i,j$ with sum $n$, the $(i,j)$-slice of the coproduct of an 
$n$-letter word $W$, written \cop{i}{j}{W}, is the Sweedler sum
$$
\sum_{(W)_{i,j}} W_{(1)}\tens W_{(2)},
$$
that is, it is the sum of tensor products of terms obtained by
decomposing the word $W$ into two subwords
of indicated lengths, both in the linear order 
induced from that on $W$, with a sign equal to that of the permutation
obtained by passing from $W$ to its reordering as the concatenation of the two subwords.
For instance, 
$$
\copm{2}{2}{abcd}= ab\tens cd - ac\tens bd + ad\tens bc + 
bc\tens ad - bd\tens ac + cd\tens ab.
$$

We define the {\it Whitney algebra of a matroid $M(S)$} as follows. 
Let $E$ be the exterior algebra freely generated by the underlying set $S$
of points of the matroid $M$, $T$ the direct sum of tensor powers of $E$,
and $W$ the quotient of $T$ by the homogeneous ideal generated by
coproduct slices of words formed from {\it dependent} sets of points.

Note that this is a straight-forward analogue of the principle applied
in the bracket algebra of a Peano space, that {\it dotting a dependent
set of points yields zero}.

This construction of a Whitney algebra
is reasonable because {\it these very identities hold in the
tensor algebra of any vector space}. Consider, for instance, a set of four
coplanar points $a,b,c,d$ in a space of rank $n$, say for large $n$.
Since $a,b,c,d$ form a dependent set, the coproduct slice 
\cop{2}{2}{abcd} displayed above will have $ij\tens kl$ -coordinate
equal to the determinant of the matrix whose rows are the coordinate
representations of the four vectors, 
calculated by Laplace expansion with respect to 
columns $i,j$ versus columns $k,l$. This determinant is zero because the 
vectors in question form a dependent set. Compare reference \cite{lecl},
where this investigation is carried to its logical conclusion.

The Whitney algebra of a matroid $M$ has a geometric product
reminiscent of the meet operation in the Grassmann-Cayley algebra,
but the product of two extensors is not equal to zero when the 
extensors fail to span the entire space. The definition is as follows,
where $A$ and $B$ have ranks $r$ and $s$, the union  $A\cup B$ has rank $t$,
and $k=r+s-t$, which would be the ``normal'' rank 
of the intersection of the flats
spanned by $A$ and $B$ if they were to form a modular pair. 
The {\it geometric
product}  is defined as either of the Sweedler sums
$$
A\diamond B = \sum_{(A)_{r-k,k}} A_{(1)}B\tens A_{(2)}= \sum_{(B)_{k,s-k}} AB_{(2)}\tens B_{(1)}
$$
This product of extensors is always non-zero. The terms in the first tensor position are
individually either equal to zero or they span the flat obtained as the span (closure) $E$ of 
the union $A\cup B$. For a represented matroid, the Grassmann coordinates of the
left-hand terms in the tensor products are equal up to an overall scalar multiple,
because $A_{(1)}B$ and $AB_{(2)}$ are non-zero if and only if 
$A_{(1)}\cup B$ and $A\cup B_{(2)}$,
respectively, are spanning sets for the flat $E$. These extensors with 
span $E$ now act like
scalars for a linear combination of \exts{k} representing the meet of $A$ and $B$,  equal to
that meet whenever $A$ and $B$ form a modular pair, and equal 
to the vector space
meet in any representation of the matroid.

The development of the Whitney algebra began with an exciting exchange of email
with \gc\ in the winter of 1995-6. On 18/11/95 he called to say that he agreed that
the ``tensor product'' approach to non-spanning syzygies is correct, so that
$$
a\tens bc - b\tens ac + c\tens ab
$$ is the zero tensor whenever $a,b,c$ are collinear points in any space, and gives a 
Hopf-algebra structure on an arbitrary matroid, potentially replacing the ``bracket ring''\cite{whiteBR}, 
which had the disadvantage of being commutative.

 Four days later he wrote:
\quo{I just read your fax, it is exactly what I was thinking. I have gone a little further 
in the formalization of the Hopf algebra of a matroid, so far everything checks beautifully. 
The philosophical meaning of all this is that every matroid has a natural coordinatization 
ring, which is the infinite product of copies of a certain quotient of the free exterior 
algebra generated by the points of the matroid (loops and links allowed, of course). This 
infinite product is endowed with a coproduct which is not quite a Hopf algebra, but a new 
object closely related to it. Roughly, it is what one obtains when one mods out all 
coproducts of minimal dependent sets, and this, remarkably, gives all the exchange 
identities. I now believe that everything that can be done with the Grassmann-Cayley 
algebra can also be done with this structure, especially meets.}

On 29/11/95: \quo{I will try to write down something tonight and send it to you by latex. 
I still think this is the best idea we have been working on in years, and all your past 
work on syzygies will fit in beautifully.}

On 20/12/95: \quo{I am working on your ideas, trying to recast them in letterplace language. 
I tried to write down something last night, but I was too tired. Things are getting quite 
rough around here.}

Then, fortunately for this subject, the weather turned bad. On 9/1/96: \quo{Thanks for the message. 
I am snowbound in Cambridge, and won't be leaving for Washington until Friday, at least,
so hope to redraft the remarks on Whitney algebra I have been collecting. \dots}

\quo{Here are some philosophical remarks. First, all of linear algebra should be done with
the Whitney algebra, so no scalars are ever mentioned. Second, there is a new theorem
to be stated and proved preliminarily, which seems to be a vast generalization of the
second fundamental theorem of invariant theory (Why, Oh why, did I not see this before?!)}

Here, \gc\ suggests a comparison between the Whitney algebra of a vector space $V$, when
viewed as a matroid, and the exterior algebra of $V$.

\quo{ I think this is the first step
towards proving the big theorem. It is already difficult, and I would appreciate your help.
The point is to prove classical determinant identities, such as Jacobi's identity, using
only Whitney algebra methods (with an eye toward their quantum generalizations!) Only by
going through the Whitney algebra proofs will we see how to carry out a quantum 
generalization of all this stuff.}

\quo{It is of the utmost importance that you familiarize yourself with the letterplace 
representation of the Whitney algebra, through the Feynman operators, and I will write 
this stuff first and send it to you.}

On 11/1/96, still snowbound in Cambridge, \gc\ composed a long text proposing two projects:
\begin{enumerate}
\item{the description of a module derived from a Whitney algebra $W(M)$,}
\item{a faithful
representation of a Whitney algebra as a quotient of a supersymmetric letter-place algebra.}
\end{enumerate}

This supersymmetric algebra representation is as follows. He uses the supersymmetric
letter-place algebra $\hbox{Super}[S^-|P^+]$ in negative letters and positive places.
\leave{the latter to represent the various positions in the tensor product}. 
A tensor product
$W_1\tens W_2\tens \dots \tens W_k$ is sent to the product
$$
(W_1|p_1^{|W_1|})(W_2|p_2^{|W_2|})\dots(W_k|p_k^{|W_k|})
$$
where the words $p_i^{(|w_i|)}$ are divided powers of positive letters, representing the different
possible positions in the tensor product. The letter-place pairs are thus anti-commutative. The 
linear extension of this map on $W(M)$ he termed the {\it Feynman entangling operator}.

Bill Schmitt joined the project in the autumn of '96. The three us met together only once.
Bill managed to solve the basic problem, showing that the Whitney algebra is precisely a {\it
lax Hopf algebra}, the quotient of the tensor algebra of a free exterior algebra by an ideal
that is not a co-ideal. The main body of this work Bill and I finished
\cite{whit} 
 in 2000, too late to share
the news with \gc. It was not until eight years later that Andrea Brini, Francesco Regonati,
Bill Schmitt and I finally established that \gc\ had been completely correct about the 
super-symmetric representation of the Whitney algebra.

\begin{figure}[h] 
     \centering
     \includegraphics{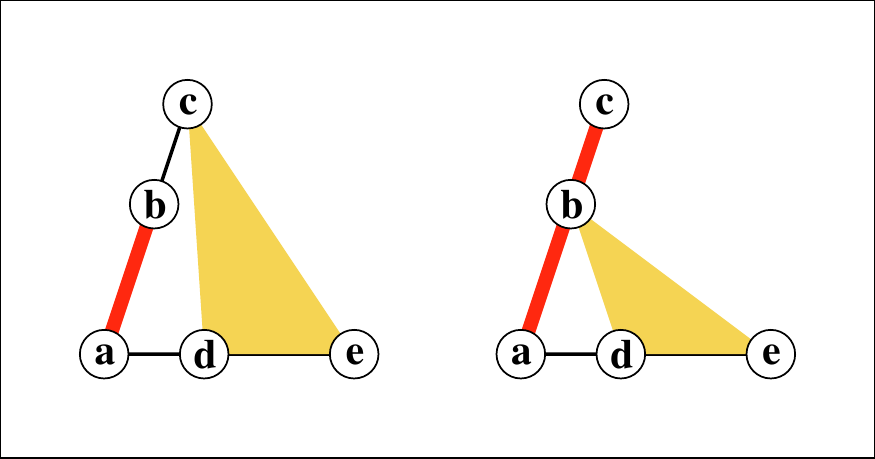}
     \caption{$ab\tens cde = ac \tens bde$.}
     \label{Fi:L}
     \end{figure}

For a quick taste of the sort of calculations one does in the Whitney algebra of a matroid,
and of the relevant geometric signification, consider a simple matroid $L$ on five points:
two coplanar lines $bc$ and $de$, meeting at a point $a$.
In Figure~\ref{Fi:L} we exhibit, for the matroid $L$, the geometric reasoning 
behind the equation
$ab \tens cde = ac \tens bde$.

First of all,
$$
ab\tens c - ac\tens b + bc\tens a = 0
$$
because $abc$ is dependent. Multiplying by $de$ in the second tensor position,
$$
ab\tens cde - ac\tens bde + bc\tens ade = 0,
$$
but the third term is zero because $ade$ is dependent in $L$, so we have
the required equality $ab \tens cde = ac \tens bde$.

This equation of tensor products expresses the simple fact 
that the ratio  $r$ of the lengths of the oriented line segments $ac$ and $ab$
is equal to the ratio of the oriented areas of the triangular regions
$cde$ and $bde$, so, in the passage from the product on the left
to that on the right, there is merely a shift of a scalar factor $r$ from the
second term of the tensor product to the first. 
{\it As tensor products they are equal.}

The same fact is verified algebraically as follows. Since $b$ is collinear
with $a$ and $c$, and $d$ is collinear with $a$ and $e$, we may write 
$b =  (1-\alpha) a + \alpha\, c$ and $d = \beta\, a + (1-\beta) e$ for 
some choice of non-zero 
scalars $\alpha,\beta$. Then 
$$\begin{matrix}
ab\tens cde = (a\join ((1-\alpha) a + \alpha\, c))\tens 
(c\join (\beta\, a + (1-\beta) e)\join e)\\
ac\tens bde = (a \join c)\tens
(((1-\alpha) a + \alpha\, c)\join (\beta\, a + (1-\beta) e)\join e)
\end{matrix}$$
both of which simplify to $-\alpha\beta\,(ac\tens ace)$

Note that this equation is independent of the dimension of the overall space
within which the triangular region $ace$ is to be found. 

\section{Geometric product}

For words $u,v\in W^1$, with $|u|=r, |v|=s$, let $k= r+s-\rho(uv)$.
The {\it geometric product} of $u$ and $v$ in $W$, written $u\diamond v$, 
is given by the expression
$$
u\diamond v = \sum_{(u)_{r-k,k}}u_{(1)}v \circ u_{((2)}
$$

For words $A,B$ and integers $r,s,k$ as above,
$$
 \sum_{(A)_{r-k,k}}A_{(1)}B \circ A_{((2)} =  \sum_{(B)_{k,s-k}}A_{(1)}B \circ A_{((2)}.
$$
So the geometric product is commutative:
$$
A\diamond B = (-1)^{(r-k)(s-k)}B\diamond A
$$

\begin{figure}[h] 
     \centering
     \includegraphics{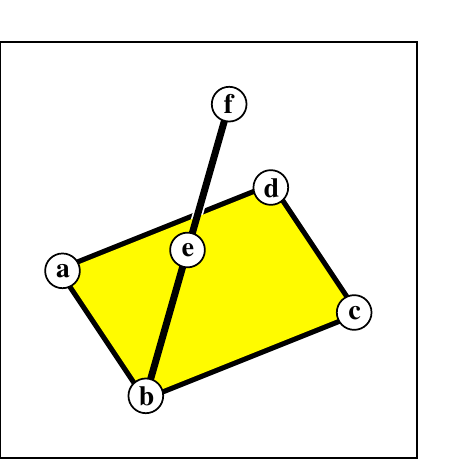}
     \caption{The geometric product of a line with a plane. }
     \label{Fi:gp1}
     \end{figure}

In Figure~\ref{Fi:gp1}, we see how the intersection point (at $b$) of the line
$ef$ with the plane $acd$ can be computed as a linear 
combination of points $e$ and $f$,
or alternately as a linear combination of points $a$, $c$ and $d$, using the two
formulations of the geometric product of $ef$ and $acd$.

The calculation is as follows:
$$
\begin{matrix}
acd\diamond ef& =& acef\circ d - adef\circ c + cdef\circ a\\
&=& acdf \circ e - acde \circ f
\end{matrix}
$$

\begin{figure}[h] 
     \centering
     \includegraphics{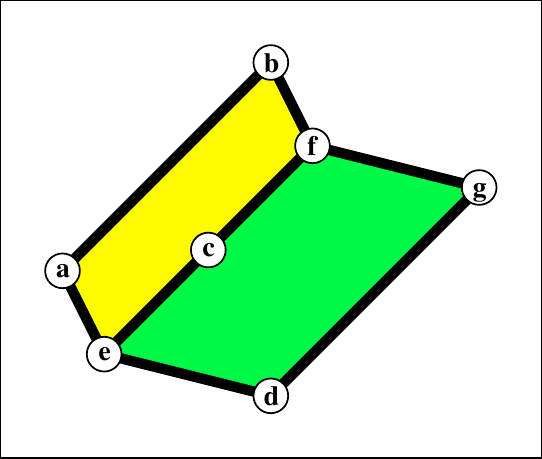}
     \caption{The geometric product of two planes. }
     \label{Fi:gp2}
     \end{figure}

Figure~\ref{Fi:gp2} shows how the line of intersection (incident with
points $a,c,f$) of planes $abc$ and $def$ can be computed as a linear
combination of lines $bc$ and $ac$, or can be obtained as a single term by the
alternate form of the geometric product.

The calculation is
$$
\begin{matrix}
abc\diamond def &=& adef\circ bc - bdef\circ ac\\
&=& abcd\circ ef
\end{matrix}$$

Compare section 56-57, pages 88-89 in Peano\cite{peanoT}.

\section{Regressive product}
When I first lectured on \gc's work on the Grassmann-Cayley algebra to a seminar
at McMaster University run by two eminent algebraists, Evelyn Nelson and Bernard Banaschewski,
they insisted that the meet operation was {\it just} the dual of Grassmann's exterior product,
and referred me to Bourbaki. I was never comfortable with this view, since I felt that
a veritable geometric product would not be restricted to the two {\it extremities} of the
lattice of subspaces, yielding non-trivial results only for independent joins and co-independent 
meets.

 Lloyd Kannenberg performed an outstanding service to mathematics when he published
his English translations \cite{grass44T,grass62T,peanoT}
 of these seldom-consulted works by Grassmann and Peano, the two
{\it Ausdehnungslehre}, published under the titles {\it A New Branch of Mathematics} (1844)
 and {\it Extension Theory} (1862) and the {\it Geometric Calculus} of Peano, 
which, as Kannenberg says,
``was published in a small print run in 1888, and has never been reissued in its entirety.''
Lloyd tells me that his translations 
of these classics were undertaken with  \gc's active encouragement.

The {\it Ausdehnungslehre} of 1844 has an entire chapter on what Grassmann calls
the regressive product. At the outset (\P125) Grassmann explains that he wants
a multiplication that will produce a non-zero value for the product of
magnitudes $A,B$ that are dependent upon one another. 
\quo{\it In order to discover this new definition we must investigate the different degrees
of dependence, since according to this new definition the product of two dependent
magnitudes can also have a nontrivial value.} (We will put all direct quotations from
the English translation in {\it italics}. We also write $\pr$ for the regressive product,
this being somewhat more visible than Grassmann's period ``.'' notation.)

To measure the different degrees of dependence, Grassmann argues that the set of points
dependent upon both $A$ and $B$ forms a {\it common system}, what we now call
a projective subspace, the intersection of the spaces spanned by $A$ and by $B$.
\quo{\it To each degree of dependence corresponds a type of multiplication: we include
all these types of multiplication under the name regressive multiplication.}
The {\it order} of the multiplication is the value chosen for the rank of the
common system.

In \P126 Grassmann studies the modular law for ranks:
$$
\ra A + \ra B = \ra C + \ra D
$$
where $\ra$ is the rank function,  $C$ is the common system 
(the lattice-theoretic meet) and $D$ 
is the nearest covering system (the lattice-theoretic join). In \P129 he 
explains the {\it meaning} of a geometric product.
\quo{ In order to bring the actual value of a real regressive product
under a simple concept we must seek, for a given product whose value
is sought, all forms into which it may be cast, without changing
its value, as a consequence of the formal multiplication principles
determined by the definition. Whatever all these forms have in common will
then represent the product under a simple concept.} So the meaning
of the regressive product is synonymous with the equivalence relation
\quo{have the same geometric product}. 

He sees that the simplest form of a product is one in {\it subordinate form},
that is when it is a {\it flag} of extensors. He thus takes the {\it value},
or meaning, of the product to be the \quo{combined observation} of the
flag of flats \quo{together with the (scalar) quantity to be distributed
on the factors.} A scalar multiple can be transferred from term to term in
a product without changing the value of the product, that is, he is introducing
a tensor product.

As a formal principle he permits the dropping of 
an additive term in a factor
if that term has a higher degree of dependence on the other factors of the 
product. For instance, in the figure of three collinear points $a,b,c$ together
with a line $de$ not coplanar with $a,b,c$, $ab\pr(ce+de) = ab\pr de$. We will
subsequently show that this product is also equal to $1\pr abde$.

\P130 gives the key definition. \quo{\it If $A$ and $B$ are the two factors of a regressive product and
the magnitude $C$ represents the common system of the two factors, 
then if $B$ is set equal to $CD$, $AD$ represents the nearest covering system
and thus the relative system as well if the product is not zero.} That is, we represent
one of the factors, $B$,
 as a product $CD$ of an extensor $C$ spanning the intersection of the subspaces spanned by
 $A$ and $B$, times an extensor $D$ that is complementary to $C$ in the 
subspace spanned by $B$. We then transfer the factor $D$ to the multiplicative term
involving $A$. The result is a flag of extensors, which Grassmann decides to write in 
decreasing order. He concludes that this flag expression is unique (as tensor product).
In \P131: \quo {\it The value of a regressive
product consists in the common and nearest covering system of the two factors, 
if the order of the factors is given, apart from a quantity to be 
distributed multiplicatively on the two systems.}

Also in \P131, he states that the regressive product of two extensors is {\it equal} 
to its associated flag representation. With $A,B$ and $B=CD$ as above, he writes
$$
A\pr B=A\pr CD=AD\pr C
$$
Perhaps even more clearly (\P132), he states that  \quo{\it we require that
two regressive products of nonzero values $A\pr CD$ and $A'\pr C'D'$ 
are equal so
long as generally the product of the outermost factors with the middle one
is equal in both expressions, 
or if they stand in reciprocal proportion, whether
the values of the orders of the corresponding factors agree or not. 
In particular,
with this definition we can bring that regressive product 
into subordinate form.}
That is, any regressive product of extensors $A\pr B$ is equal to a flag product,
say $E\pr C$, where $E=AD$ in the previous calculation.

\begin{figure}[h] 
     \centering
     \includegraphics{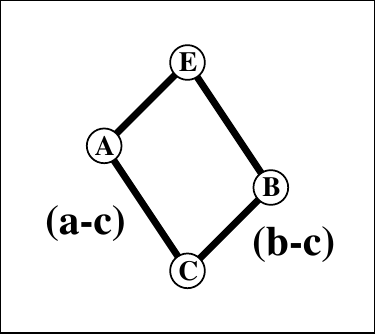}
     \caption{Values of rank in a modular lattice. }
     \label{Fi:mod}
     \end{figure}

He also figures out the sign law for commutativity of the regressive product.
If $A\pr B=E\pr C$, and the right hand side is a flag, where the ranks of the extensors
$A,B,C,E$ are $a,b,c,e$, respectively, then the sign for the exchange of order
of the product is the parity of the product $(a-c)(b-c)$, the product of the 
{\it supplementary numbers} (See Figure~\ref{Fi:mod}). It's fascinating that he managed to get this right!

So we have (skew) commutativity for the regressive product. {\it How about associativity?}
Grassmann calls associativity {\it the law of combination.} 
This is the extraordinary part of the story.

He has a well defined product (in fact, one for every degree of dependence), so, in principle,
any multiple product is well defined (and even the correct degrees for each successive
product of two factors), but he soon recognizes, to his dismay, that the regressive
product is {\bf not}  generally associative.

\begin{figure}[h] 
     \centering
     \includegraphics{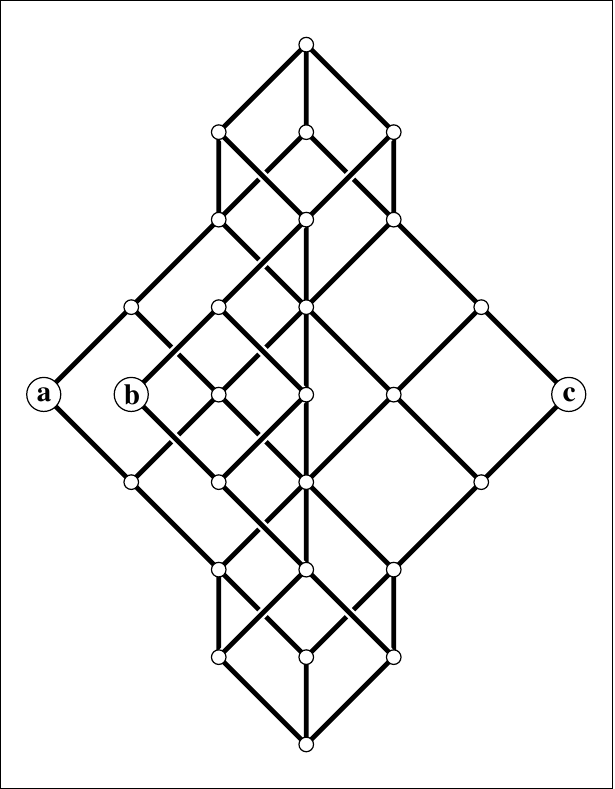}
     \caption{The free modular lattice on $3$ generators. }
     \label{Fi:ded}
     \end{figure}

The best way to see the non-associativity is by reference to the free modular lattice
on three generators, as first found by Dedekind, in Figure~\ref{Fi:ded}. 
(The elements $a,b,c$ are indicated with lower-case letters, but they could 
be extensors of any rank. I'll also write flags in increasing order.)

In any lattice, for any pair of elements in the order $x<z$, and
for any element $y$, 
$$
x\join(y\meet z)\le (x\join y)\meet z
$$ 
A lattice is {\it modular} if and only if, under these same conditions, 
equality holds:
$$
x\join (y\meet z) =  (x\join y)\meet z
$$
The lattice of subspaces of a vector space is modular
(and complemented), so any calculus of linear
pieces of space will use the logic of modular lattices.

\begin{figure}[h] 
     \centering
     \includegraphics{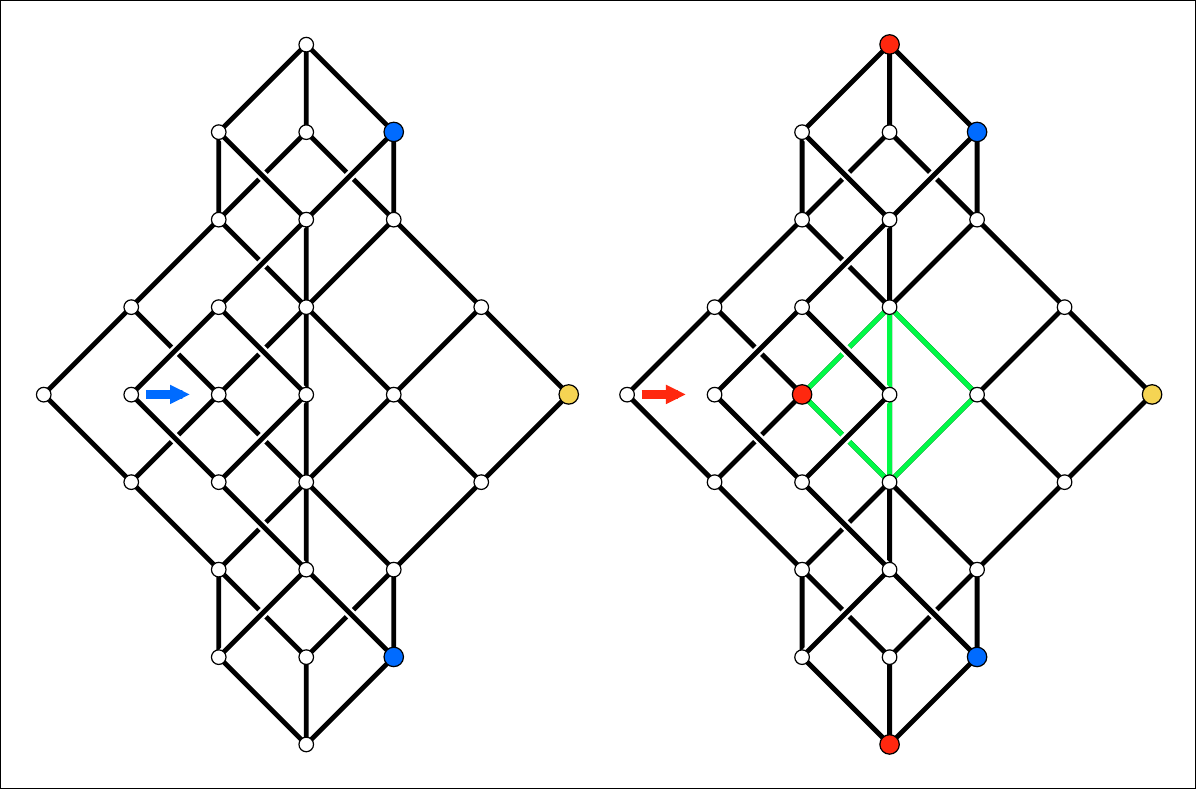}
     \caption{The regressive product of $a$ with $b\diamond c$. }
     \label{Fi:ded-bca}
     \end{figure}

On the left of Figure~\ref{Fi:ded-bca}, we indicate\col{ in blue} the passage from 
$b\pr c$ to $b\meet c\pr b\join c$.
On the right of Figure~\ref{Fi:ded-bca}, we show\col{ in red} the passage from
$b\meet c\pr b\join c$ to 
$$
a\pr ((b\meet c)\pr (b\join c)) = 
(a\meet b\meet c)\pr((a\join(b\meet c))\meet(b\join c))\pr(a\join b\join c)
$$
 A value of the triple product will always land at one of the
three central elements of the inner\colAlt{ green}{ gray} sublattice, but the exact
position will depend on which factor entered last into the combined product.
{\it The result carries a trace of the order in which the factors were combined!}

Try this with the simple case of three collinear points $a,b,c$. Then
$$
a\pr(b\pr c) = a\pr (1\pr bc) = 1\pr a\pr bc, 
$$
but
$$
(a\pr b)\pr c = (1\pr ab)\pr c = 1\pr c\pr ab. 
$$

Parentheses are not necessary in the final flags, because flag products are
associative.
All products in which $a$ enters last will be the same, up to sign:
$a\pr(b\pr c) = (b\pr c)\pr a = - (c\pr b)\pr a$.

However, the regressive product of two flag products is well defined.
Grassmann showed this to be true, {\it an incredible feat, given the tools he
had at hand!} In fact, he proved that the extensors of one flag can progressively
multiplied into another flag, and that the result is independent not only
of the order of these individual multiplications, but also independent
of which flag is multiplied into which! In \P 136 Grassmann says \quo{\it
Instead of multiplying by
a product of mutually incident factors one can progressively multiply by
the individual factors, and indeed in any order.
}
 
Garrett Birkhoff, in the first edition of his {\it Lattice Theory}, proved that
the free modular lattice generated by two finite chains is
a finite {\it distributive} lattice. This changes the whole game.

\begin{figure}[h] 
     \centering
     \includegraphics{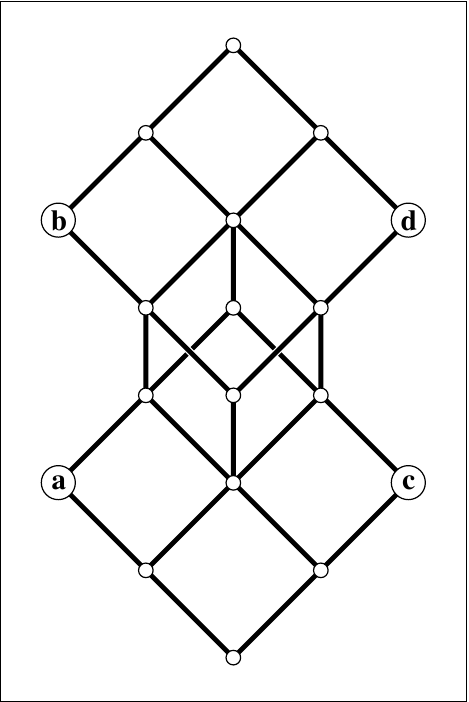}
     \caption{The free modular lattice generated by two $2$-chains. }
     \label{Fi:twoCh}
     \end{figure}

\begin{figure}[h] 
     \centering
     \includegraphics{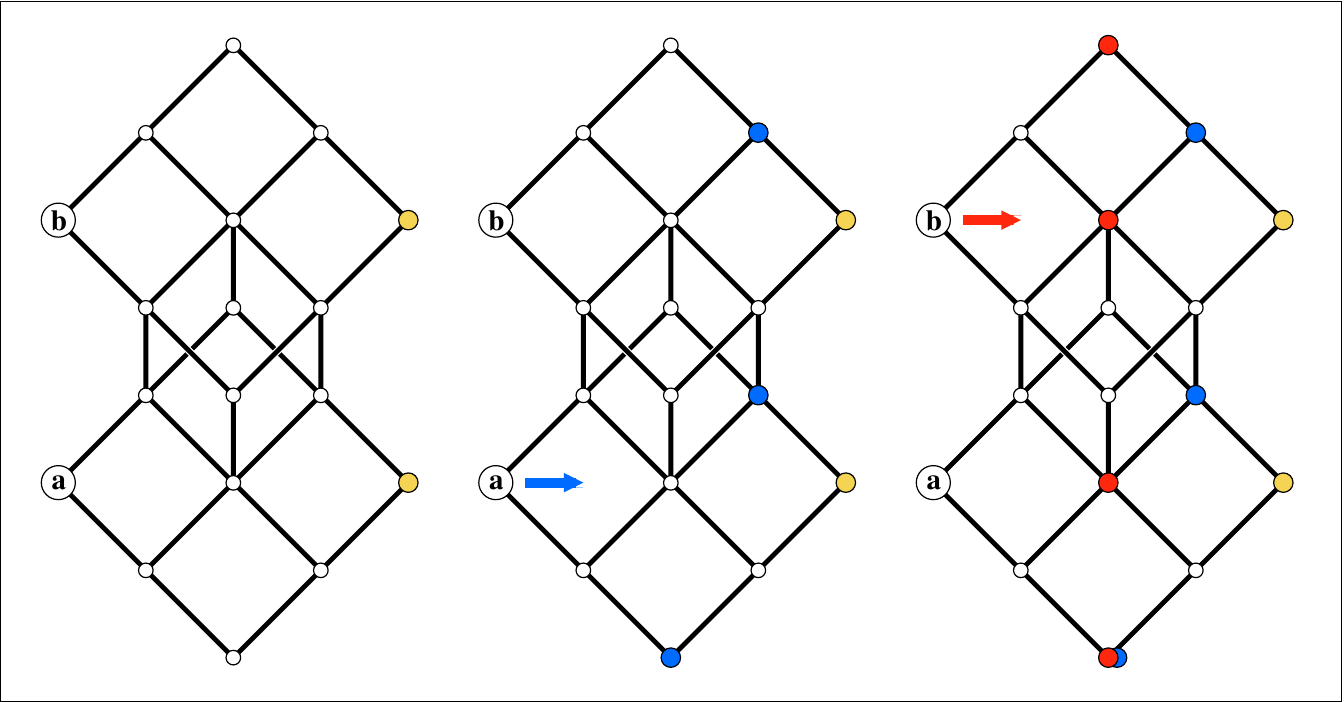}
     \caption{The product  $b\to(a\to (c<d))$}
     \label{Fi:twoCh-ab}
     \end{figure}

\begin{figure}[h] 
     \centering
     \includegraphics{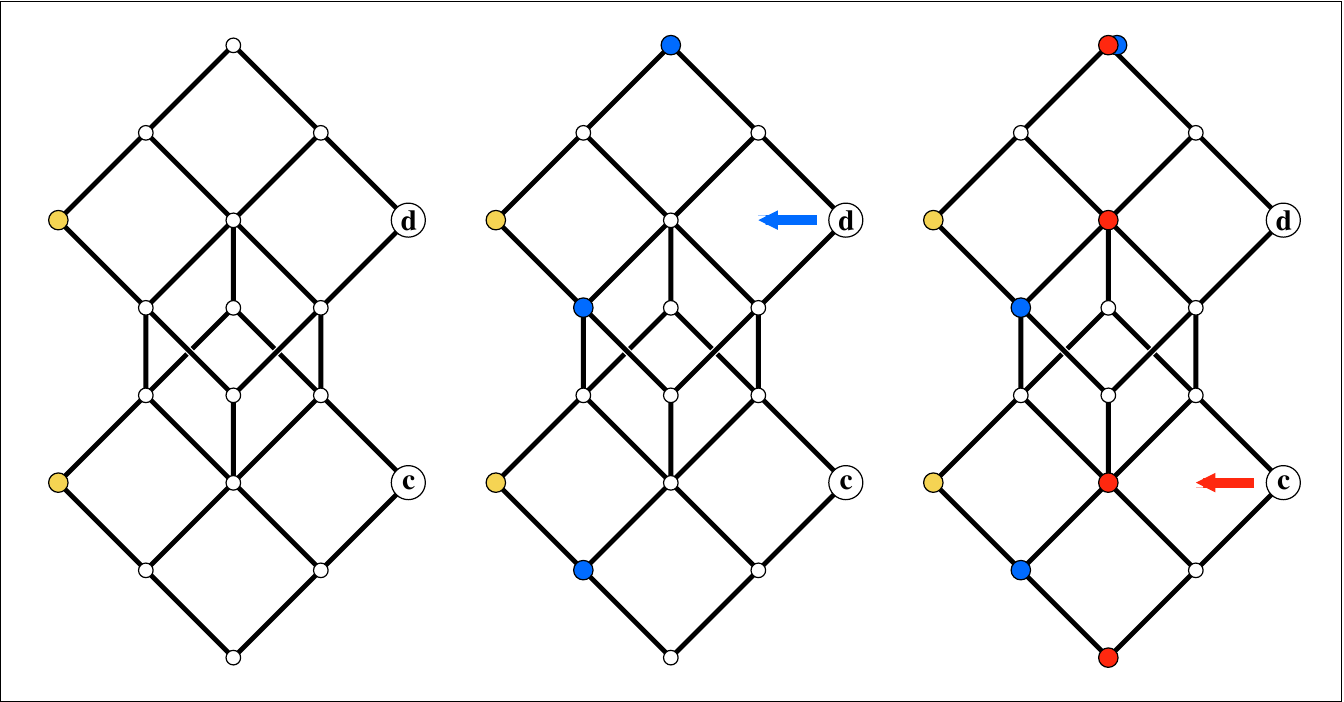}
     \caption{The product  $c\to(d\to (a<b))$. }
     \label{Fi:twoCh-dc}
     \end{figure}

In Figure~\ref{Fi:twoCh} we show the free modular lattice generated by two $2$-chains.
In the center of Figure~\ref{Fi:twoCh-ab} we show the result of multiplying the extensor
$a$ into the flag $c\pr d$, then, on the right, the result of multiplying $b$ into that 
result. We end up on what might well be termed the {\it backbone} of the 
distributive lattice. In Figure~\ref{Fi:twoCh-dc} we show what happens when
the flag $c\pr d$ is multiplied into the flag $a\pr b$, but starting with the 
top element $d$, instead. The result is the same. 

By the end of his chapter on the regressive product, Grassmann seems rather disheartened.
He admits clearly in \P139 that \quo{\it the multiplicative law of combination \dots is not
generally valid for the mixed product of three factors.}

He includes a footnote
to say that {\it cases can be found in which our law still finds its application}
via the results available for the product of an extensor by a flag of extensors, but
concludes, with a certain degree of disillusion, that \quo{these cases
are so isolated, and the constraints under which they occur so contrived, 
that no scientific advantage results from their enumeration.}

Then, having investigated duality and having proven associativity for joins of
independent extensors and meets of co-independent extensors, he concludes with the note: 
\quo{\it the theoretical 
presentation of this part of extension theory now appears as completed, 
excepting consideration of the types of multiplication for which the law of 
combination is no longer valid.} 
He adds the footnote: 
\quo{\it How to treat such products, which to be sure have manifold applications, 
I have sought to indicate at the conclusion of this work.}

Clarification of these questions involving the regressive product of two flags has
been joint work with Andrea Brini, Francesco Regonati, and Bill Schmitt, last year
in Bologna.

When this work is combined with the extraordinary synthesis
\cite{btProd}  of Clifford algebra and
Grassmann-Cayley algebra, all made super-symmetric, already achieved by Andrea Brini, 
Paolo Bravi, and Francesco Regonati, here present, you have finally the 
makings of a banquet that
can truly be termed {\it geometric algebra}.

\section{Higher Order Syzygies}
Before closing, we should take a quick look at the {\it higher order syzygies}
that \gc\ mentioned in his email messages concerning the new Whitney algebra.

Given a configuration $C$ of $n$ points in projective space
of rank $k$, certain subsets
$A\subseteq C$ will be dependent. The minimal dependent sets 
(the {\it circuits} of the corresponding matroid $M(C)$) will have dependencies
that are uniquely defined up to an overall scalar multiple. They thus form,
in themselves, a configuration of rank $n-k$ of projective points in a space 
of rank $n$. We call this the first {\it derived configuration, $\derc 1$}, 
and denote the associated matroid  $\derm 1$, 
the {\it derived matroid}. 

In the same way, 
the circuits of $\derc 1$ form a new projective configuration, which we 
denote $\derc 2$, with matroid $\derm 2$, and so on. Thus, any matroid
represented as a configuration in projective space automatically
acquires an infinite sequence of derived matroids. In classical
terminology, $\derc k$ is the configuration of {\it $k^{th}$-order 
syzygies} of $C$.

\begin{figure}[h] 
     \centering
     \includegraphics{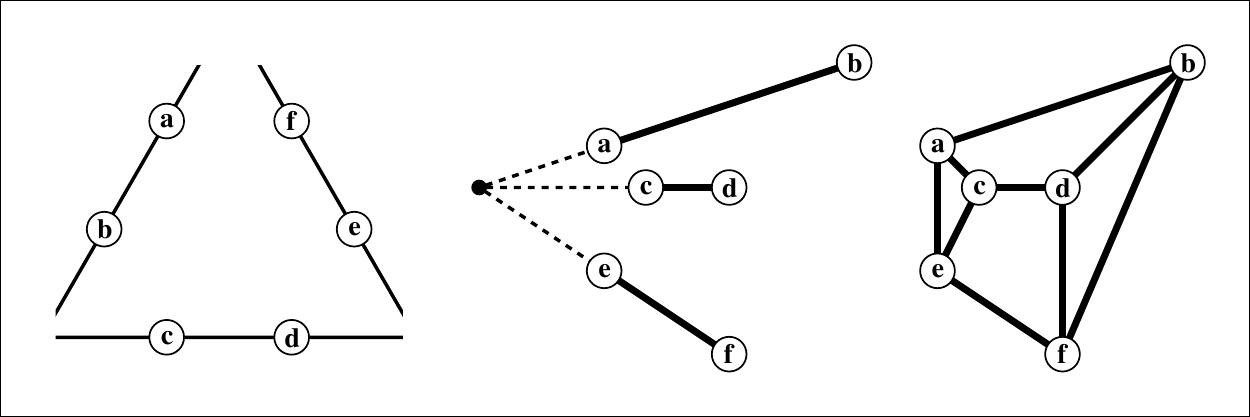}
     \caption{Representations of the uniform matroid $U_{3,6}$. }
     \label{Fi:U36}
     \end{figure}

This derived information is {\it not}, however, fully determined by
the matroid itself. The simplest example of interest is given by the uniform
matroid $U_{3,6}$ of six points $\{a,b,\dots, f\}$ in general position 
in the projective plane. In Figure~\ref{Fi:U36}, 
we show two representations of
the matroid $U_{3,6}$ in the plane. In the example on the left, the three 
lines $ab, cd, ef$ do not meet, and the circuits $abcd, abef, cdef$ are 
independent, spanning the space of first order syzygies among the six points.
On the right, those three lines meet, the circuits $abcd, abef, cdef$ are
dependent (rank $2$ in the derived matroid). Those three circuits act
as linear constraints on lifting of the figure of six points into $3$-space.
A height vector $(h(a), h(b), \dots, h(f))$ is orthogonal to the vector
$([bcd], [acd], [abd], [bcd], 0, 0)$ if and only if the four lifted points
$a'=(a_1,a_2,a_3,h(a)), \dots, d'=(d_1,d_2,d_3,h(d))$ are coplanar, 
if and only if {\it the dependency is preserved in the lifting}. If the three
circuits are of rank 3, there will be $6-3$ choices for the lift of the
six points, and they must remain coplanar. The three circuits are of rank $2$
if and only if the plane figure has a polyhedral lifting, to form a true
three-dimensional triangular pyramid. And this happens if and only if the 
three lines $ab, cd, ef$ meet at a point, namely, the projection of the
point of intersection of the three distinct planes 
$a'b'c'd', a'b'e'f', c'd'e'f'$.

\begin{figure}[h] 
     \centering
     \includegraphics{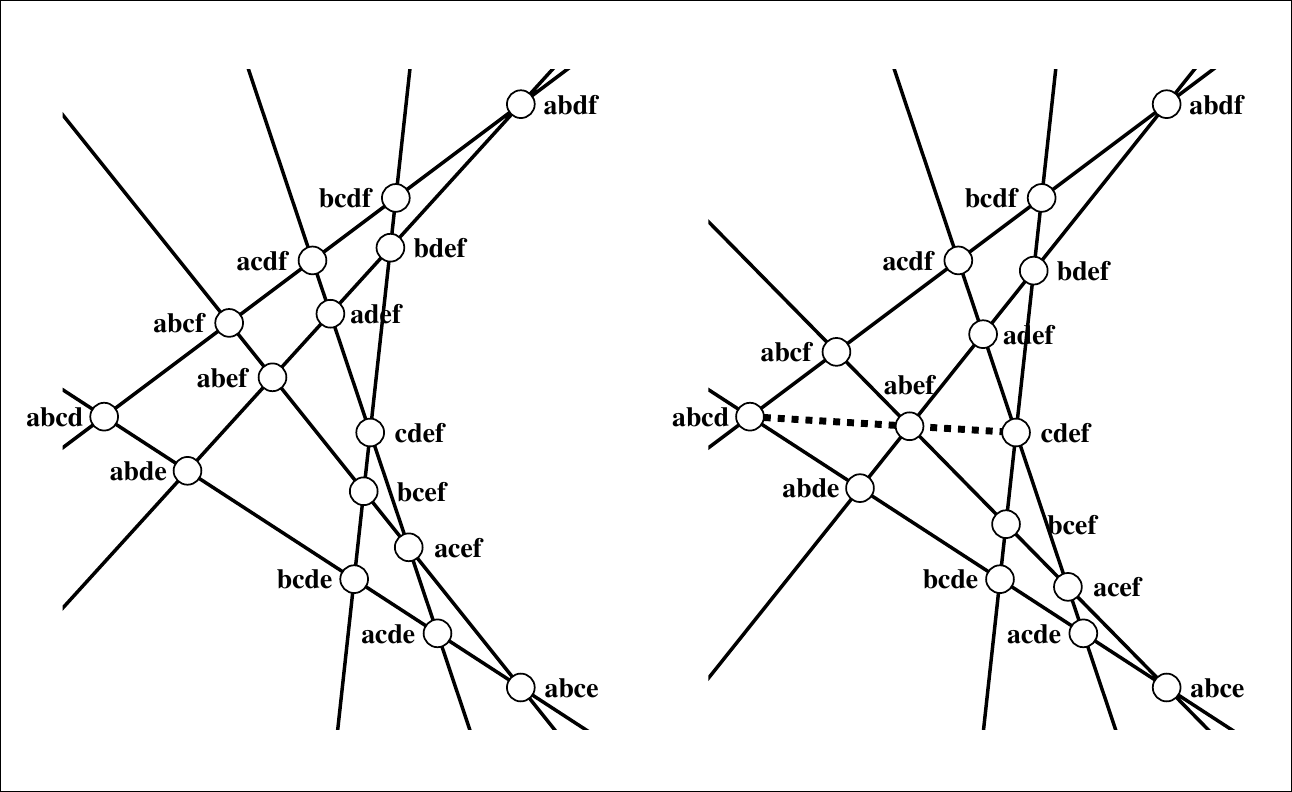}
     \caption{First order syzygies for the uniform matroid $U_{3,6}$. }
     \label{Fi:sixD}
     \end{figure}

For the figure on the left, the derived configuration of first order
syzygies is of rank $3$ and consists simply of the intersection of 
6 lines in general position in the plane. The five circuits formed from
any five of the six points have rank $2$, and are thus collinear, as in
Figure~\ref{Fi:sixD}. 
In the special position, when the lines $ab, cd, ef$ are
concurrent, the circuits $abcd, abef, cdef$ become collinear, as on the
right in Figure~\ref{Fi:sixD}.

The geometric algebra can also be put to service to provide coefficients for
higher-order syzygies, in much the same way that generic first-order syzygies
are obtained as coproducts of dependent sets. To see how this is done, it will suffice
to take another look at the present example. Check the matrix of coefficients
of these first-order syzygies:
$$\begin{matrix}
&a&b&c&d&e&f\\
(abcd)&abc&-acd&abd&-abc&0&0\\
(abef)&bef&-aef&0&0&abf&-abe\\
(cdef)&0&0&def&-cef&cdf&-cde
\end{matrix}$$
This row space $S$ is orthogonal to the space of {\it linear functions} on
the set $S=\{a,b,c,d,e,f\}$, a space $T$ spanned by the coordinate projection
functions
$$\begin{matrix}
&a&b&c&d&e&f\\
(1)&a_1&b_1&c_1&d_1&e_1&f_1\\
(2)&a_2&b_2&c_2&d_2&e_2&f_2\\
(3)&a_3&b_3&c_3&d_3&e_3&f_3
\end{matrix}$$
The Grassmann coordinates of these two vector subspaces differ 
from one another
(up to an overall ``scalar'' multiple) by {\it complementation of places} 
and a sign of that complementation. 
This algebraic operation is called the  {\it Hodge star} operator.
Concretely, where $\sgn$ is the sign of the permutation merging two
words into one, in a given linear order,
$$\begin{matrix}
S_{abc}=\sgn(def,abc) T_{def}=-T_{def}\\
S_{abd}=\sgn(cef,abd) T_{cef}=T_{cef}\\
\dots\\
S_{def}=\sgn(abc,def) T_{abc}=T_{abc}\\
\end{matrix}
$$
The Grassmann coordinate $T_{xyz}$ for any letters $x,y,z$ is just the
\ext3\ $xyz$.
The ``scalar'' in question, which \gc\ and I called {\it the resolving bracket}
\cite{resBra} is obtainable by calculating any $3\times 3$ minor
of the matrix for $S$, and dividing by the \ext3\ obtained from the complementary
set of columns in the matrix for $T$. If we do this in columns $abc$ for $S$,
we find determinant $(-bcd\tens aef + acd\tens bef) def$, which, when divided by
$-def$ yields $bcd\tens aef - acd\tens bef$. It helps to recognize that this expression
can be obtained by joining the meet of $ab\meet cd = a\tens bcd - b\tens acd$ with $ef$, 
so this resolving bracket is equal to zero if and only if the meet of $ab$ and $cd$ is
on the line $ef$, that is, if and only if the three lines $ab, cd, ef$ are concurrent.
This is the explicit synthetic condition under which the three first-order syzygies
$abcd, abef, cdef$ form a dependent set in the first derived configuration.

Much work remains in order to develop an adequate set of tools for dealing with
higher order syzygies in general. The concept of resolving bracket is but a first
step. \gc\ and I spent many hours discussing these higher order syzygies, usually
on the white-boards in his apartment in Boston, later in Cambridge, in his office,
or in more exotic places such as Strasbourg or Firenze, during mathematical gatherings.
 I think he enjoyed these discussions, in the period 1985-95,
difficult as it was for him to force me to express
my ideas clearly. The only major breakthrough was in \gc's
very fruitful collaboration with David Anick \cite{ar,arC}, 
where they found a resolution of the
bracket ring of a free exterior algebra, bases for syzygies of all orders being
represented by families of totally non-standard tableaux. In this way, you have only to deal
with syzygies having {\it single bracket coefficients}.

\section{Balls in Boxes}
As a closing thought, I would like to express my conviction 
that \gc\ was rightfully
fascinated by probabilistic questions arising from quantum theory, but somehow
never really got a proper hold on the basic issues, despite having approached
them from all quarters: via general combinatorial theory, 
{\it esp\`eces de structure}, 
supersymmetric algebra, umbral calculus, probability theory, and \dots philosophy.

Let me suggest that it is high time we reread what he has written
here and there 
on {\it balls and boxes}, as in title of today's memorial assembly, 
for hints he may generously have left us.

As he wrote in the introduction to {\it The Power of Positive Thinking} (with Wendy Chan),

\narrower{
\quo{\it The realization that the classical system of 
Cartesian coordinates can and 
should be enriched by the simultaneous use of 
two kinds of coordinates, some of which 
commute and some of which anticommute, 
has been slow in coming; its roots, like the 
roots of other overhaulings of our prejudices 
about space, go back to physics, to 
the mysterious duality that is found in 
particle physics between identical particles 
that obey or do not obey the Pauli exclusion principle.}}



\begin{thebibliography}{99}

\bibitem{ar}
D.~Anick, G.-C.~Rota,
\emph{Higher-order Syzygies for the Bracket Ring and for the Ring of Coordinates of the Grassmannian,}
Proc. Nat. Acad. of Sci. \textbf{88} (1991), 8087-8090.

\bibitem{bbr}
M.~Barnabei, A.~Brini, G.C.~Rota,
\emph{On the Exterior Calculus of Invariant Theory,}
J. of Algebra \textbf{96} (1985), 120-160.

\bibitem{bb}
P.~Bravi, A.~Brini, 
\emph{Remarks on Invariant Geometric Calculus, Cayley-Grassmann Algebras and Geometric Clifford Algebras,}
in H.~Crapo, D.~Senato, \emph{Algebraic Combinatorics and Computer Science, 
A Tribute to Gian-Carlo Rota,}
Springer 2001, pp 129-150

\bibitem{bht}
A.~Brini, R.~Q.~Huang, A.~G.~B.~Teolis, 
\emph{The Umbral Symbolic method for Supersymmetric Tensors,}
Adv. Math., \textbf{96} (1992), 123-193.

\bibitem{brt}
A.~Brini, F.~Regonati, A.~G.~B.~Teolis,
\emph{Grassmann Geometric Calculus, Invariant Theory and Superalgebras,}
in H.~Crapo, D.~Senato, \emph{Algebraic Combinatorics and Computer Science, 
A Tribute to Gian-Carlo Rota,}
Springer 2001, pp 151-196.

\bibitem{btProd}
A.~Brini, A.~Teolis,
\emph{Grassmann's Progressive and Regressive Products and GC Coalgebras, }
in G.~Schubring, ed., \emph{Hermann G\"unther Grassmann (1809-1877), Visionary Mathematician, Scientist and Neohumanist Scholar,} Kluwer (1996), 231-242.

\bibitem{chan3-6}
W.~Chan, 
\emph{Classification of Trivectors in $6-D$ Space},
in \emph{Mathematical Essays in Honor of Gian-Carlo Rota,}
B.~E.~Sagan and R.~P.~Stanley, ed., Birkhäuser 1998, pp~63-110.

\bibitem{crs}
W.~Chan, G.-C.~Rota, J.~Stein,
\emph{The Power of Positive Thinking,} in
\emph{Proceedings of the Curaçao Conference: Invariant Theory in Discrete and Computational Geometry,} 1994, Kluwer, 1995.

\bibitem{resBra}
H.~Crapo, G.-C.~Rota,
\emph{The Resolving Bracket}
in \emph{Proceedings of the Curaçao Conference: Invariant Theory in Discrete and Computational Geometry}, 1994, Kluwer, 1995.

\bibitem{arC}
H.~Crapo
\emph{On the Anick-Rota Representation of the Bracket Ring of the Grassmannian,}
Advances in Math., \textbf{99} (1993), 97-123.

\bibitem{whit}
H.~Crapo, W.~Schmitt,
\emph{The Whitney Algebra of a Matroid,}
J.of Comb. Theory (A), 

\bibitem{drs}
P.~Doubilet, G.-C.~Rota, J.~Stein, 
\emph{On the Foundations of Combinatorial Geometry: IX, 
Combinatorial Methods in Invariant Theory}. Studies in Applied mathematics, No 3, vol LIII, Sept 1974, pp 185-215,.

\bibitem{grass44}
Hermann Grassmann,
\emph{Die lineale Ausdehnungslehre, ein neuer Zweig der Mathematik,}
Verlag von Otto Wigand, Leipzig, 1844. 

\bibitem{grass44T}
Hermann Grassmann,
\emph{A New Branch of Mathematics: The Ausdehnungslehre of 1844, and Other Works,}
translated by Lloyd C. Kannenberg,
Open Court, 1995.

\bibitem{grass62}
Hermann Grassmann,
\emph{Die Ausdehnungslehre, Vollst\"andig und in strenger Form,}
Verlag von Th. Cgr. Fr. Enslin (Adolph Enslin), Berlin, 1862.

\bibitem{grass62T}
Hermann Grassmann,
\emph{Extension Theory}
translated by Lloyd C. Kannenberg,
Amer. Math. Soc. 2000.

\bibitem{grs}
F.~Grosshans, G.-C.~Rota, J.~Stein,
\emph{Invariant Theory and Supersymmetric Algebras,}
Conference Board of the Mathematical Sciences, No 69,
Amer. Math. Soc., 

\bibitem{hrs}
R.~Huang, G.-C.~Rota, J.~Stein,
\emph{Supersymmetric Bracket Algebra and Invariant Theory,}
Centro Matematico V. Volterra, Università Degli Studi di Roma II, 1990.

\bibitem{lecl}
B.~Leclerc,
\emph{On Identities Satisfied by Minors of a Matrix,}
Advances in Math., \textbf{100} (1993), 101-132.

\bibitem{peanoT}
Giuseppe Peano
\emph{Geometric Calculus, according to the Ausdehnungslehre of H. Grassmann,}
translated by Lloyd C. Kannenberg,
Birkhäuser, 2000.

\bibitem{whiteBR}
N.~White,
\emph{The Bracket Ring of a Combinatorial Geometry, I and II,}
Trans. Amer. Math. Soc. \textbf{202} (1975a) 79-95, \textbf{214} (1975b) 233-48.

\bibitem{whiteCGA}
N.~White,
\emph{A Tutorial on Grassmann-Cayley Algebra, }
Proceedings of the Curaçao Conference: Invariant Theory in Discrete and Computational Geometry, 1994, Kluwer, 1995.

\end{thebibliography}
\end{document}